\documentclass[11pt]{article}

\usepackage[T2A]{fontenc}
\usepackage[utf8]{inputenc}
\usepackage[english,russian]{babel}

\usepackage{amsmath,amssymb,amsthm,amsfonts,mathrsfs}
\usepackage{graphicx}

\usepackage[margin=1in]{geometry}
\usepackage{microtype}
\usepackage[hidelinks]{hyperref}

\usepackage{placeins}   
\usepackage{flafter}    
\usepackage{needspace}  
\usepackage{float}      
\title{\bfseries {\LARGE On a Configuration with Circle-Conic Tangency and Sharygin Points}}
\author{
Petr Kim\thanks{Liceum Vorobievy Gory, Moscow, Russia, Independent University of Moscow}
\and
Georgii Makoian\thanks{Letovo School, Moscow, Russia, email: makoyan2008@gmail.com}
}\date{} 

\begin{document}
\maketitle

\small \bf Abstract. \it This work studies circle-geometry methods through their application to a main theorem about circles tangent twice to a conic. The authors investigate the Sharygin point---a point lying in the pencil of two non-intersecting circles---and explore its properties. These properties are applied to solve several olympiad problems, such as problems from MGO 2024 and the Croatian IMO selection. The paper also presents a simplified version of the main theorem and gives two different proofs: one using Sharygin points and another using Lobachevsky (hyperbolic) geometry. The article explores relation between Lorenz transformations of Minkowski space-time and a certain transformation of circles in hyperbolic geometry. The paper demonstrates the effectiveness of combining classical planimetry with ideas of non-Euclidean geometry for solving difficult problems involving circle tangencies.
\vspace{.3cm}
\normalsize
\Needspace{6\baselineskip}\nopagebreak
\begin{center}
\large \bf 0. Main Results
\end{center}
\vspace{.3cm}
\rm The following theorem showcases up-to-date methods of classical geometry:\\
\\
\bf Main Theorem\\
\it Let $\omega, \omega_1$ be intersecting circles passing through points $A,B$ and $C,D$, respectively. The points $A,B,C,D$ are concyclic. Let $\gamma$ be a conic passing through $A,B,C,D$ whose center lies inside both $\omega$ and $\omega_1$. If there exists a circle tangent to $AB$, $CD$, $\omega$, and $\omega_1$, touching one of $\omega,\omega_1$ internally and the other externally, then there exists a circle tangent to $\gamma$ at two points and tangent to $\omega$ and $\omega_1$. \\
\\
\begin{figure}[H]
\centering
\includegraphics[width=0.6\linewidth]{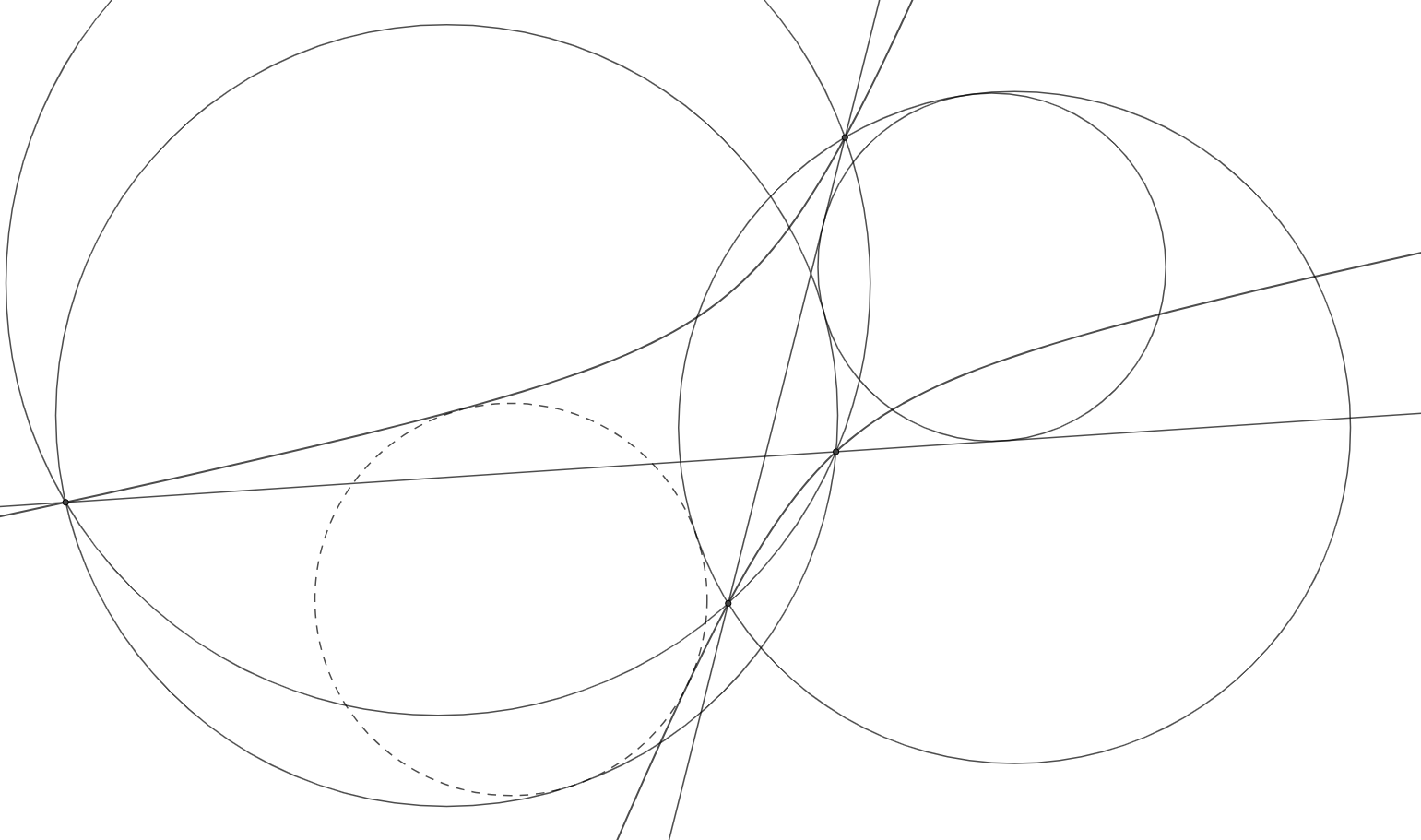}
\end{figure}
\noindent\rm We also consider related facts:\\
\\
\bf Simplified Main Theorem\\
\it Let $\omega, \omega_1$ be intersecting circles passing through points $A,B$ and $C,D$, respectively. Then the points $A,B,C,D$ are concyclic. If there exists a circle tangent to $AB$, $CD$, $\omega$, and $\omega_1$, touching one of $\omega,\omega_1$ internally and the other externally, then for any circle $\Gamma$ tangent to $\omega$ and $\omega_1$ with one internal and the other external tangency, there exists a circle $\gamma$ tangent to $AB$ and $CD$ and lying in the pencil $(ABCD), \Gamma$.\\
\Needspace{6\baselineskip}\nopagebreak
\begin{center}
\large \bf 1. Proof of the Weak Main Theorem via Properties of the Sharygin Point of Two Circles
\end{center}
\vspace{.2cm}
\rm The use of Sharygin points will be demonstrated on the following olympiad problems:\\
\\
\bf 1. MGO 2024, Problem 4, author: Grigory Zabaznov\\
\it In triangle $ABC$, points $D$, $E$ on sides $AB$, $AC$ are such that quadrilateral $BCED$ is cyclic. Let $K$ be the intersection point of $BE$ and $CD$. Points $L$, $M$ on $BE$, $CD$ respectively are such that $A$, $L$, $M$ lie on the line symmetric to $AK$ with respect to the bisector of $\angle BAC$. Prove that $(KLM)$ is tangent to $DE$ if and only if it is tangent to $BC$.\\
\\
\bf 2. Croatian IMO Selection 2016, Problem 7, author unknown\\
\it In triangle $ABC$ a point $S$ satisfies $\frac{AS+BS}{AB} = \frac{BS+CS}{BC} = \frac{AS+CS}{CA}$. Let $A_1$, $B_1$, $C_1$ be the second intersections of $AS$, $BS$, $CS$ with $(ABC)$. Prove that the incircles of $\triangle ABC$ and $\triangle A_1B_1C_1$ coincide.\\
\\
\bf 3. Weak Main Theorem:\\
\it Let $\omega, \omega_1$ be circles passing through vertices $A,B$ and $C,D$ respectively of a cyclic quadrilateral $ABCD$. Prove the equivalence of the following statements: a) there exists a circle tangent to $AB$, $CD$, $\omega$, $\omega_1$; b) there exists a circle concentric with $(ABCD)$ tangent to $\omega$, $\omega_1$, with one external and the other internal tangency; c) if $S$ is an intersection point of $\omega$ and $\omega_1$, then in the pencil containing $S$ and $(ABCD)$ there exists a circle tangent to $AB$ and $CD$.\\
\\
\begin{figure}[H]
\centering
\includegraphics[width=0.6\linewidth]{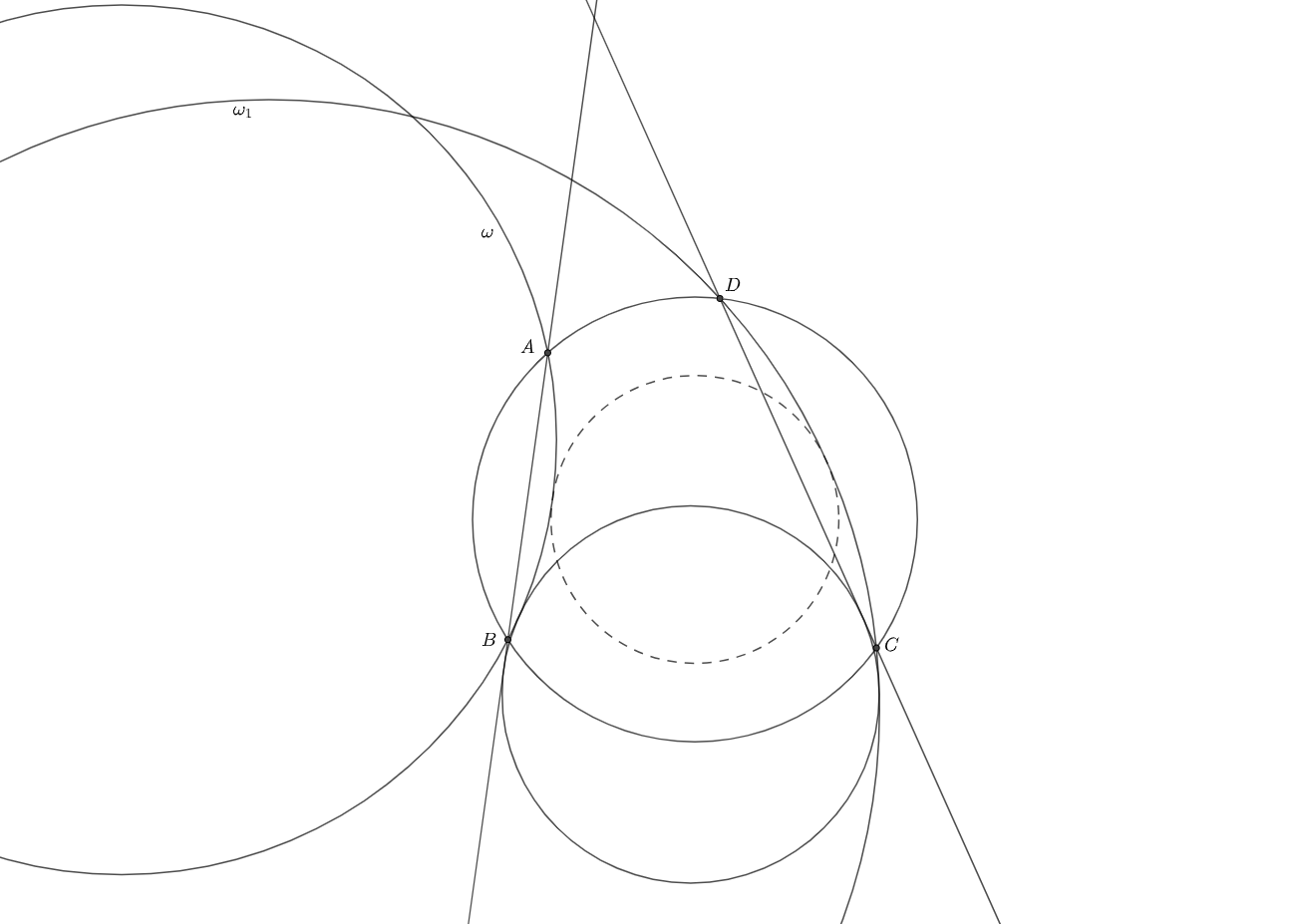}
\end{figure}
\FloatBarrier
\newpage
\noindent \rm A slightly modified version of this theorem was proposed at the Southern Mathematical Tournament 2023.\\
\\
\bf Definition 1.1\\
\it We call a \emph{Sharygin point} of two non-intersecting circles $\omega_1 , \omega_2$ any point lying in the pencil containing $\omega_1 , \omega_2$. (From the equation of a coaxial pencil it is clear that a pair of circles has two such points.) \\
\\
\bf Lemma 1.1\\
\it Let $\omega_1$ and $ \omega_2$ be two non-intersecting circles, and let
$ S$ and $ S'$ be their Sharygin points.  
Then the polars of $S$ with respect to $\omega_1$ and $\omega_2$ coincide, and this common polar is the line through $S'$ perpendicular to $O_1O_2$.\\
\\
\bf Proof\\
\rm Clearly, $S \in O_1O_2$. The perpendicular bisector of $SS'$ is the radical axis of any pair of circles in the pencil. In particular, it is the radical axis of the point $S$ and circle $\omega_1$, hence its image under the homothety with center $S$ and ratio $2$ is the polar of $S$ with respect to $\omega_1$. Therefore the line through $S'$ perpendicular to $O_1O_2$ is the polar of $S$ with respect to $\omega_1$. Similarly, this line is also the polar of $S$ with respect to $\omega_2$.\\
\\
\bf Lemma 1.2\\
\it Let $\omega_1$ and $ \omega_2$ be two non-intersecting circles, and let
$ S$ be their Sharygin point.  
Then the images of $\omega_1$ and $ \omega_2$ under inversion with center $S$ and radius $p$ are concentric circles.\\
\\
\bf Proof\\
\it
Let $O_1'$ and $O_2'$ be the centers of the images of $\omega_1$ and  $ \omega_2$ under inversion. As is well known,
\[
\overrightarrow{SO_i'} \;=\; \frac{\rho^{2}}{\,SO_i^{2}-r_i^{2}\,}\,\overrightarrow{SO_i},
\qquad i\in\{1,2\}.
\]
We also observe that 
\[SO_1 - 
\frac{r_1^{2}}{SO_1} = SO_1 - O_1S' = SS' = SO_2 - O_2S' = SO_2 - \frac{r_2^{2}}{SO_2}.
\]
Hence
\[
\overrightarrow{SO_1'} = \overrightarrow{SO_2'}  ,
\]
which yields the claim.\\
\\
\\
\rm From these lemmas we obtain one of the principal properties of the Sharygin point:\\
\\
\bf Property 1\\
\it Let $S$ be a Sharygin point of $\omega_1,\omega_2$, and let a chord $AB$ (or its extension)
of $\omega_1$ intersect $\omega_2$ at $C,D$.
Let $A',B'$ be the second intersections of $SA,SB$ with $\omega_1$, and
$C',D'$ the second intersections of $SC,SD$ with $\omega_2$.
Then $A',B',C',D'$ are collinear.\\
\\
\bf Proof\\
\rm From the polar property of the Sharygin point it follows that under a projective transformation sending $S$ to the common center of the images of $\omega_1,\omega_2$, the circle $\omega_2$ becomes a conic with the same center. The statement is then immediate by symmetry. \\
\\
\bf Corollary 1\\
\rm Under the assumptions of Property 1, if $AB$ is tangent to $\omega_2$,
then the line $A'B'$ is tangent to $\omega_2$.\\
\\
\bf Property 2\\
\it Let $S$ be a Sharygin point of $\omega_1,\omega_2$, and let a chord $AB$ (or its extension)
of $\omega_1$ be tangent to $\omega_2$ at $X$.
Then $X$ lies on the bisector of $\angle ASB$.\\
\\
\bf Proof\\
\it We use a standard lemma on pencils of circles:
for points $U,V$ there exists a circle through $U,V$ belonging to the pencil
$(S,\omega_2)$ if and only if
\[
\frac{Pow_S(U)}{Pow_{\omega_2}(U)}=\frac{Pow_S(V)}{Pow_{\omega_2}(V)}.
\]
Apply this to $U=A$, $V=B$, taking into account that $AB$ is tangent to $\omega_2$ at $X$:
$Pow_{\omega_2}(A)=XA^2$, $Pow_{\omega_2}(B)=XB^2$.
We obtain
\[
\frac{SA^2}{XA^2}=\frac{SB^2}{XB^2}\;\;\Longrightarrow\;\;
\frac{SA}{XA}=\frac{SB}{XB}.
\]
In triangle $SAB$ this is equivalent to $SX$ being an angle bisector.\\
\\
\bf Property 3\\
\it Let $S$ be a Sharygin point of $\omega_1,\omega_2$, and let a chord $AB$ (or its extension)
of $\omega_1$ be tangent to $\omega_2$.
Then the quantity $\displaystyle \frac{SA+SB}{AB}$ is independent of the choice of $A,B$.
Moreover, if for points $C,D\in\omega_1$ one has
$\displaystyle \frac{SA+SB}{AB}=\frac{SC+SD}{CD}$, then $CD$ is tangent to $\omega_2$.\\
\\
\bf Proof\\
\it From Property 2 we have $\dfrac{SA}{XA}=\dfrac{SB}{XB}$,
where $X$ is the point of tangency of $AB$ with $\omega_2$.
\vspace{0.1cm}
Since $X$ lies on the line $AB$, we have $AB=XA+XB$, hence
\[
\frac{SA+SB}{AB}
=\frac{SA}{XA}=\frac{SB}{XB},
\]
which yields the independence from the choice of $A,B$.
The converse is proved similarly.\\
\\
\bf Property 4\\
\it Let $S$ be a Sharygin point of $\omega_1,\omega_2$, and let a chord $AB$ (or its extension)
of $\omega_1$ meet $\omega_2$ at $C,D$.
Then $\angle ASC=\angle BSD$ (and therefore $\angle ASD=\angle BSC$).\\
\\
\bf Proof\\
\it For variety (it is clear this property can be proved analogously to Property 2) we give a proof via Desargues's Involution Theorem:\\
\it Consider the involution on line $AB$ which swaps points $C$ and $D$, and also swaps any two intersection points of $AB$ with a circle from the pencil of $\omega_1$ and $\omega_2$. Let $T$ be the intersection of $AB$ with the radical axis of $\omega_1$ and $\omega_2$. Then $|TS|^2 = TA\cdot TB$, hence $T$ is the foot of the tangent to the circumconic $(SAB)$. Under the involution above, $T$ maps to the point at infinity. Note that reflection in the bisector of $\angle ASB$ sends line $SA$ to $SB$, and line $ST$ to a line parallel to $AB$. Since an involution is determined by two pairs of points, the claim follows.\\
\begin{figure}[H]
\centering
\includegraphics[width=0.7\linewidth]{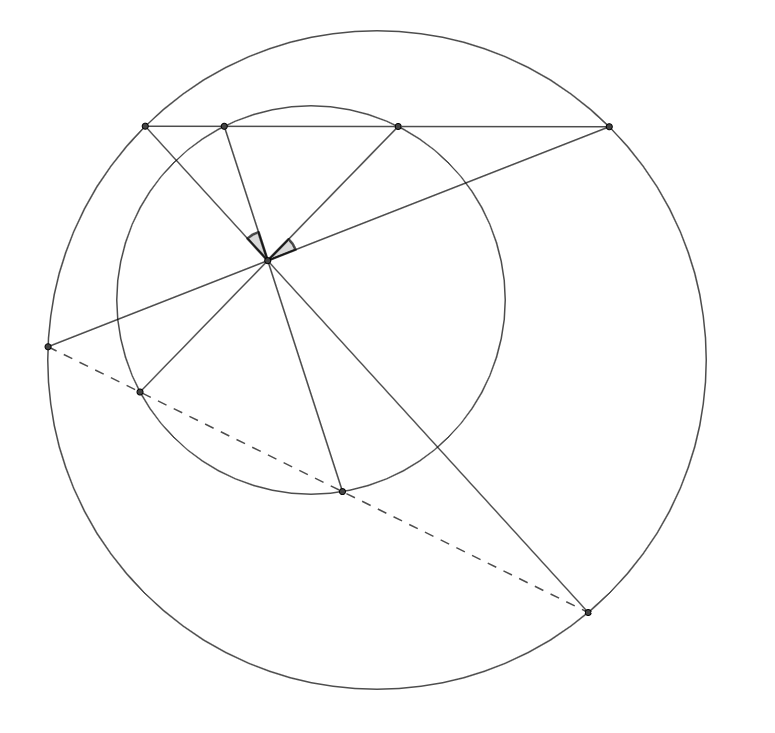}
\end{figure}
\vspace{.1cm}
\noindent\bf Remark\\
\rm There are exactly two projective homologies with center $S$ and axis $l = \operatorname{pol}_{\omega_1}(S) = \operatorname{pol}_{\omega_2}(S)$ that send $\omega_1$ to $\omega_2$ (that is, projective transformations fixing all lines through $S$ and all points on $l$ and mapping $\omega_1$ to $\omega_2$), and they are very similar.\\
\\
\bf Proof\\
\it Consider a projective transformation $\Phi$ of the plane such that $\Phi(l)=\ell_\infty$. Since polarity is preserved under projective transformations, $\Phi(S)$ is the common center of the images of $\omega_1$ and $\omega_2$. Let $\Gamma_1,\Gamma_2$ be the images of $\omega_1$ and $\omega_2$, respectively. Let $H$ be a homothety with center $\Phi(S)$ sending $\Gamma_1$ to $\Gamma_2$ (there are exactly two: with positive and with negative ratio).\\
Then consider the plane transformation $T := \Phi^{-1} \circ H \circ \Phi$. It is projective, and it maps $\omega_1$ to $\omega_2$. Since $H$ does not affect $\ell_\infty$, the line $l$ is fixed pointwise, and every line through $\Phi(S)$ is preserved by $H$; hence $T$ is a homology with center $S$ and axis $l$ taking $\omega_1$ to $\omega_2$. Since $H$ can be one of two homotheties, we have obtained two such homologies.\\
A homology is determined by its axis, center, and a pair “point–image”. Take any point $A \in \omega_1$ such that line $AS$ meets $\omega_2$ in two points $C$ and $D$ (clearly $C$ and $D$ do not lie on $l$). Since $AS$ must be sent to itself, $A$ has only two possible images, $C$ or $D$, hence exactly two such homologies exist, as required.\\
\\
\bf On the name\\
\rm It is easy to see that the 255th triangle center of $ABC$ is the Sharygin point of the incircle and one of the excircles. In English sources the 255th center is sometimes called the Sharygin point, so we found this name appropriate.\\
\\
\rm We now return to the problems stated at the beginning of the section:\\
\\
\bf 1. First proof. \\
\rm It is not hard to see that the polars of point $A$ with respect to $(BCED)$ and $(KLM)$ coincide. Hence $A$ is a Sharygin point of $(BCED)$ and $(KLM)$. The claim then follows immediately from Property $1'$.\\
\\
\bf 1. Second proof. \\
\it This proof is unrelated to the topic of the paper, but is included due to the exceptional “self-similarity” idea that appears in it.\\
\rm Consider the homothety centered at $A$ and the reflection in the angle bisector of $\angle BAC$ which sends $B \rightarrow F$ and $C \rightarrow E$. Then $K$ maps to $K_1$, $N$ to $N_1$, and $M$ to $M_1$. It is known that $(K_{1}N_{1}M_{1})$ is tangent to $EF$. Thus, if we prove that the radical axis of $(K_{1}N_{1}M_{1})$ and $(KMN)$ is $EF$, the problem will be clear. Note that $EK_1 \parallel CK$ and $FK_1 \parallel BK$. Also, $K_1$ obviously lies on line $AMN$. A simple angle chase shows that $(KMN)$ is tangent to $AK$. Then, by Fuss's lemma, $K,M,K_1,N_1$ lie on one circle. Hence the radical center of $(K_{1}N_{1}M_{1})$, $(KMN)$, $(KMK_{1},N_{1})$ is $F$. Therefore $F$ lies on the radical axis of $(K_{1}N_{1}M_{1})$ and $(KMN)$. Similarly, $E$ lies on the radical axis of $(K_{1}N_{1}M_{1})$ and $(KMN)$. This completes the proof.\\
\begin{figure}[H]
\centering
\includegraphics[width=0.8\linewidth]{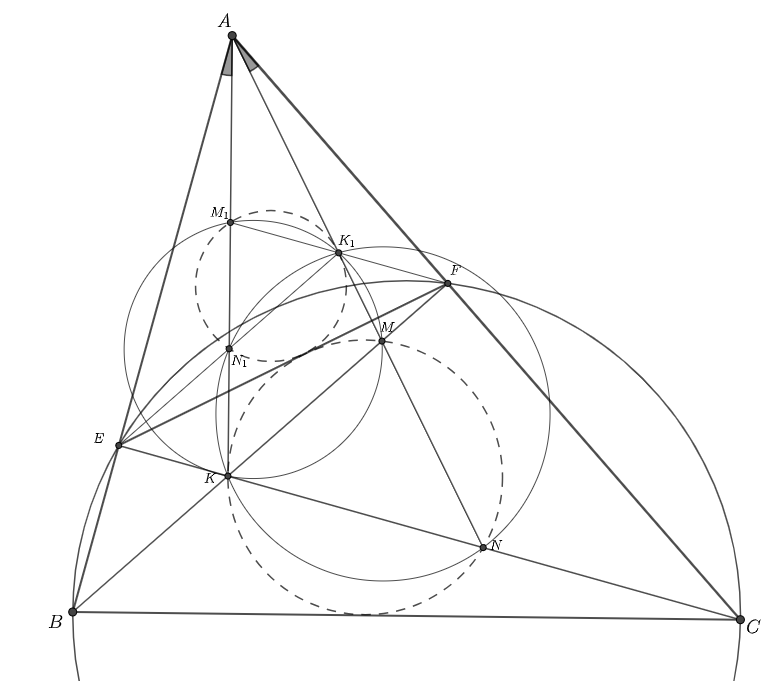}
\end{figure}
\vspace{.1cm}
\noindent\bf 2. Proof\\
\rm From Property 3 it is easy to see that $P$ is a Sharygin point of the circumcircle $(ABC)$ and the incircle of $\triangle ABC$. The statement follows from Property $1'$.\\
\\
\bf 3. Proof
\\
\rm We prove the implication from b) to c). By Property 3, $\frac{AS+BS}{AB} = \frac{CS+DS}{CD}$. Perform inversion about point $S$ with radius $R$. Then $A$ goes to $A'$, $B$ to $B'$, $C$ to $C'$, $D$ to $D'$. The circle $(ABCD)$ maps to $(A'B'C'D')$. Note that $\triangle SAB \sim \triangle SB'A'$ and $\triangle SCD \sim \triangle SD'C'$. It follows that $\frac{SA'+SB'}{A'B'}=\frac{SC'+SD'}{C'D'}$. Hence, by Property 3, there exists a circle tangent to $A'B'$ and $C'D'$ lying in the same pencil with $(ABCD)$ and $S$. Under the inverse inversion it maps to a circle tangent to $(SAB)$ and $(SCD)$, concentric with $(ABCD)$.\\
\\
\rm The implication from c) to b) becomes the previous implication under inversion in $S$. \\
\\
\rm It remains to prove the equivalence of these two statements with a).\\

\rm First we prove b) $\Rightarrow$ a). Note that $E$ is the radical center of $\omega_1$, $\omega_2$. Hence, if $S_1$ is the second intersection point of $\omega_1$ and $\omega_2$, then $S_1,S,E$ are collinear. Also, from c) we have $\frac{SA+SB}{AB} = \frac{SC+SD}{CD}$. Let $T_1$ be the foot of the angle bisector of triangle $ASB$, and $T_2$ the foot of the bisector of $CSD$. By Property $2'$, there exists a circle $\Gamma$ tangent to $AB$ and $CD$ at $T_1$, $T_2$ and lying in the pencil containing $S$ and $(ABCD)$. Let $W_1$, $W_2$ be the midpoints of arcs $AB$ and $CD$ of the circles $(SAB)$ and $(SCD)$ not containing $S$. Note that triangles $SW_{1}A$ and $AT_{1}W_{1}$ are similar. It follows that $\frac{SA+SB}{AB} = \frac{SA}{AT_1} = \frac{SW_1}{AW_1} = \frac{SW_1}{T_{1}W_{1}} = \sqrt{\frac{AW_1}{T_{1}W_{1}}}$. Similarly, $\frac{SC+SD}{CD} = \sqrt{\frac{SW_2}{T_{2}W_{2}}}$. Hence $\frac{SW_1}{T_{1}W_{1}} = \frac{SW_2}{T_{2}W_{2}}$. Therefore $\frac{ST_1}{SW_1} = \frac{ST_2}{SW_2}$, so $T_{1}T_{2} \hspace{.1cm} || \hspace{.1cm} W_{1}W_{2}$. Let $V_1$, $V_2$ be points on $AB$ and $CD$ respectively such that $EV_{1} = EV_{2}$, and let the circle tangent to $AB$ and $CD$ at $V_1$, $V_2$ touch the arc $(DSC)$ at $P$. Apply Desargues's theorem to $\triangle T_{1}W_{1}V_{1}$ and $\triangle T_{2}W_{2}V_{2}$. Note that $V_{1}V_{2} \hspace{.1cm} || \hspace{.1cm} T_{1}T_{2} \hspace{.1cm} || \hspace{.1cm} W_{1}W_{2}$. Hence, if $Y$ is the intersection of $V_{1}W_{1}$ and $V_{2}W_{2}$, then $S,E,Y$ are collinear. By Archimedes's lemma, $P$ lies on $V_{2}W_{2}$. Let $Q$ be the second intersection of $V_{1}W_{1}$ with $(SAW_{1}B)$. Since $Y$ lies on the radical axis of $(SAB)$ and $(SCD)$, points $W_{1},W_{2},P,Q$ lie on one circle. By Fuss's lemma, since $V_{1}V_{2} \parallel W_{1}W_{2}$, point $Q$ lies on the circle $(V_{1}V_{2}P)$. By Archimedes's lemma, since this circle is tangent to $AB$ at $V_1$ and passes through $Q$, it is tangent to $(SAB)$. This completes the proof.\\
\\
\rm Now prove a) $\Rightarrow$ b) by contradiction. Consider the circle concentric with $(ABCD)$ tangent to $(SAB)$. Suppose it is not tangent to $(SCD)$. Draw a circle through $C$, $D$ tangent to it. Then by the previous argument there exists a circle tangent to $(SAB)$, $AB$, $CD$, and to it. Contradiction.\\
\\
\begin{figure}[H]
\centering
\includegraphics[width=1\linewidth]{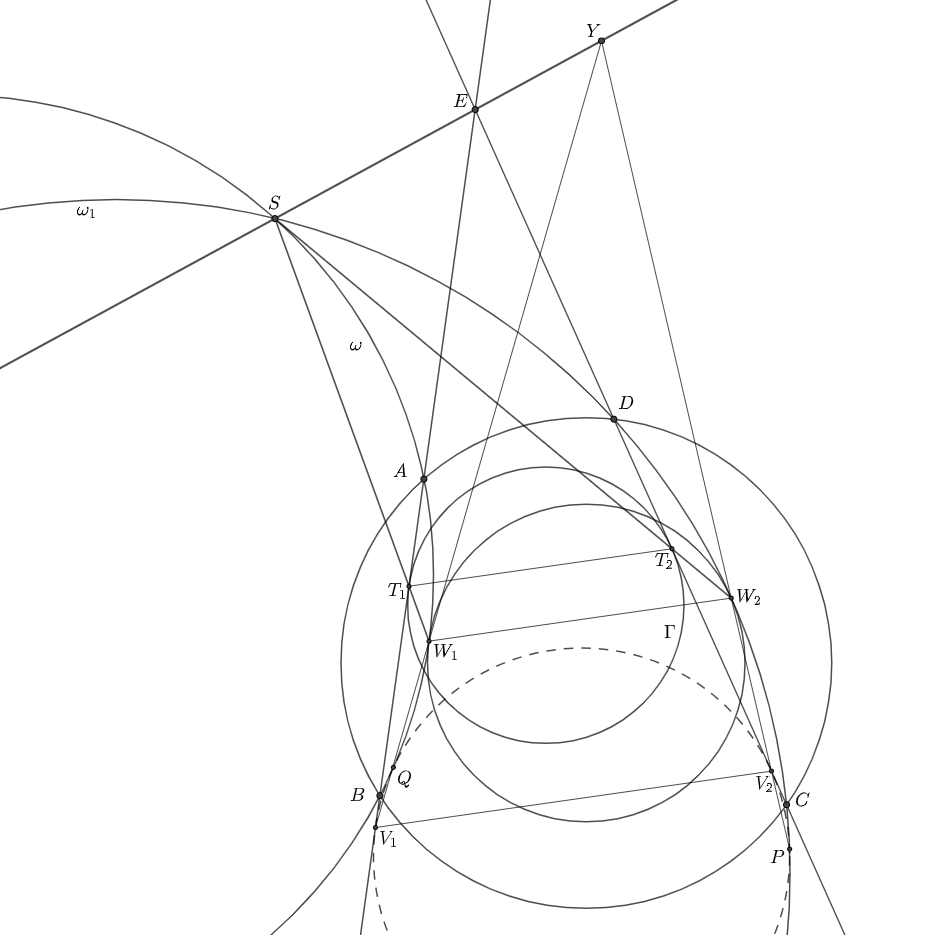}
\end{figure}
\noindent\begin{center}
\large \bf 2.1 Pseudo-Euclidean Geometry and Its Application to Euclidean Problems\\
\end{center}
\vspace{.2cm}
\it All statements in this section were proved in [1] or [2].\\
\\
\bf Definition 2.1.1\\
\it
We call a \emph{cycle} an oriented circle (clockwise or counterclockwise), and an \emph{axis} an oriented line. Two cycles are called \emph{tangent} if their circles are tangent and the tangency respects orientation. Tangency of an axis and a cycle is defined analogously. Two axes cannot be tangent. A point and an axis are called \emph{incident} if the axis passes through the point. Similarly for a point and a cycle.\\
\\
\begin{figure}[H]
\centering
\includegraphics[width=1\linewidth]{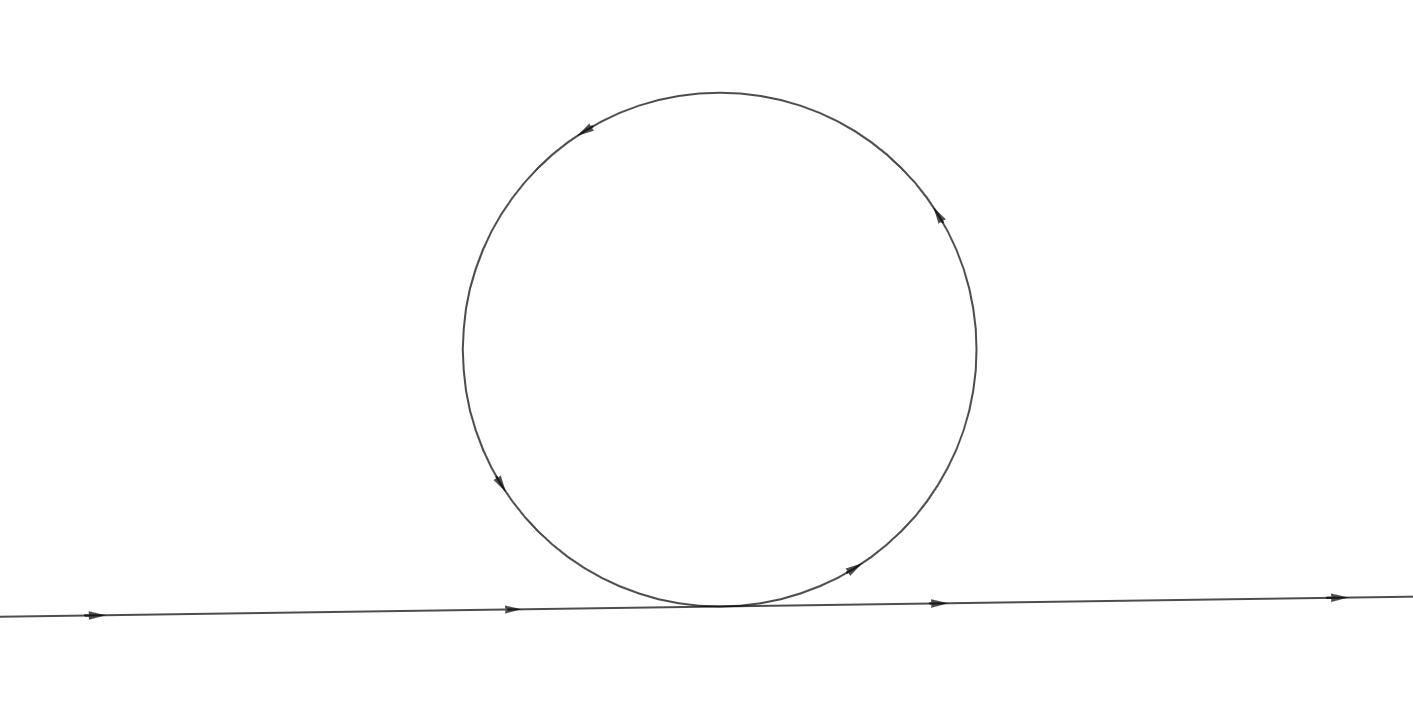}
\end{figure}
\noindent\bf Definition 2.1.2\\
\it
Define the radius of a cycle as the radius of its circle if the cycle is oriented clockwise, and as minus the radius if it is oriented counterclockwise. A point is regarded as a cycle of zero radius.\\
\\
\bf Definition 2.1.3\\
\it
A \emph{pseudo-Euclidean space} (in this paper we consider only the pseudo-Euclidean space of signature $(2,1)$) is the usual three-dimensional space $\mathbb{R}^3$ endowed with the function on pairs of points $A, B \in \mathbb{R} ^ 3$: $q(A,B) = \left( x_a - x_b \right)^2 + \left( y_a - y_b \right)^2 - \left( z_a - z_b \right)^2 $, where $\left( x_a, y_a, z_a\right)$ are the coordinates of $A$ and $\left( x_b, y_b, z_b\right)$ are the coordinates of $B$.
\\
\\
Define a bijection $\sigma$ from the set of all cycles to $\mathbb{R}^3$: to a cycle $\omega$ of radius $r$ with Cartesian center $(x,y)$ assign the point $(x,y,r)$. For brevity we will write $q(\omega, \gamma)$, where $\omega, \gamma$ are cycles, for the function $q$ of the images of the cycles in Minkowski space.
\\
Note that $q$ of tangent cycles equals $0$. Also note that if two cycles $\omega_1, \omega_2$ have a common tangent axis and the distance between its points of tangency with the cycles is $l$, then $q(\omega_1, \omega_2) = l^2$.\\
\\
\bf Definition 2.1.4\\
\it An \emph{inflation} by radius $\rho$ is the replacement of every cycle by the concentric cycle with radius increased by $\rho$, and every axis by a parallel axis at distance $\rho$ to its left. Points are regarded as cycles of radius $0$ under this transformation.\\
\\
\rm Under $\sigma$ an inflation by radius $\rho$ becomes the translation by vector $(0,0,\rho)$ in space.\\
\rm The main property of inflation is that it sends tangent objects to tangent objects.
\\
\Needspace{6\baselineskip}\nopagebreak
\begin{center}
\large \bf 2.2 Proof of the Simplified Main Theorem Using Lobachevsky Geometry
\end{center}
\it A brief introduction to hyperbolic geometry and its application to Euclidean problems can be found in [3]. We use the terminology of that paper.\\
\it Comment: the argument proposed in the previous section can be substantially simplified using the upper-half-plane Poincaré model with absolute $SE$, but the authors wished to present a proof without non-Euclidean techniques.\\
\\
\rm Call a \emph{generalized cycle} in hyperbolic geometry any oriented line, oriented circle, oriented equidistant curve (hypercycle), oriented horocycle, or a point. On the set of generalized cycles one can define an analogue of inflation, called \emph{hyperbolic inflation}:\\
\rm By analogy with Euclidean geometry, define the radius of an oriented circle as the radius of the corresponding circle if it is oriented clockwise, and as minus the radius if it is oriented counterclockwise. The clockwise direction can be fixed, for example, by using the Poincaré disk model in Euclidean geometry.\\
\rm Replace each oriented circle by a circle whose radius is increased by $\rho$, and replace each oriented equidistant curve or line by an oriented equidistant curve or line of the same direction lying in its pencil and at distance $\rho$ from it. Two equidistant curves/lines of one pencil have the same direction if, in the fixed Poincaré disk model, at their two common points the direction into/out of the disk is the same for the pair. Images of horocycles are defined analogously, except that two horocycles are called equally oriented if they are tangent as oriented circles in the Poincaré model. A point is again regarded as a circle of radius $0$.\\
\rm This transformation sends tangent generalized cycles to tangent generalized cycles (the verification is essentially identical to the Euclidean case and reduces to checking trivial cases).\\ 
\\
\rm To prove the simplified main theorem one may use the Poincaré disk model with absolute $(ABCD)$, but first we must define an “extension” of this transformation to all Euclidean cycles, axes, and points. \\
\rm Extend the transformation to the entire Euclidean plane as follows: define it on the complement of the disk by applying inversion to the transformation defined above on the interior of the disk. The action of inversion on the direction defined on a generalized cycle is arranged as follows: to each direction one associates a cyclic or linear order of points, which inversion maps to the order of points on the image of the corresponding object. Any Euclidean cycle or axis either lies entirely inside or entirely outside the disk, in which case its image is well defined and is a cycle or axis, or it lies partly inside and partly outside the disk. In that case one can observe that the images of its parts together form a Euclidean cycle. \\
\rm Proof: Let $\gamma$ be a cycle intersecting the absolute $\Omega$ at points $A,B$. Then the image of the arc of $\gamma$ lying outside the disk bounded by $\Omega$ is an equidistant curve symmetric to the arc of $\gamma$ inside the disk with respect to the line through $A$ and $B$, with opposite orientation relative to the arc inside the disk. These two arcs are matched by inversion with respect to the circle through $A$ and $B$ orthogonal to $\Omega$, since the composition of inversions in orthogonal circles sends any circle of their pencil to itself. We call two oriented equidistant curves of one pencil \emph{conjugate} if they are symmetric with respect to the line in their pencil and have opposite orientations. The two oriented arcs corresponding to the arc of $\gamma$ inside the disk and to the inverse image of the arc of $\gamma$ outside $\Omega$ form a pair of conjugate equidistant curves. Under hyperbolic inflation conjugate oriented equidistant curves map to conjugate ones. It follows that under hyperbolic inflation the image of any Euclidean cycle, axis, or point is an axis, cycle, or point.\\
\rm Under hyperbolic inflation, tangent or incident objects map to tangent or incident objects. One also notes that each Euclidean cycle maps to a cycle whose circle lies in the pencil with it and with $\Omega$.\\
\rm Now apply hyperbolic inflation to prove the simplified main theorem: \\
\rm First assign orientations to all objects. Choose an arbitrary orientation on the circle $\alpha$ tangent to $\omega_1$, $\omega$, $AB$, $CD$. The orientation of $\alpha$ uniquely determines orientations on $\omega_1$, $\omega$, $AB$, $CD$. The orientations on $\omega_1$, $\omega$ determine an orientation on $\Gamma$. Consider the hyperbolic inflation sending $\alpha$ to the oppositely oriented circle to $\alpha$. Then $\omega$ maps to $AB$ with opposite orientation to that first assigned to it. $\omega_1$ maps to $AB$ with opposite orientation to that first assigned to it. Under this transformation, $\Gamma$ maps to a cycle whose circle lies in the pencil with $\Gamma$ and $(ABCD)$ and is tangent to $AB$ and $CD$. The existence of this circle is exactly what had to be proved. \\
\\

\rm The Sharygin point also has a simple interpretation. The Sharygin point of $(ABCD)$ and $\Gamma$ is the center of the circle $\Gamma$ in the Poincaré model. From this, Property 1 is proved by central symmetry.\\
\\
\Needspace{6\baselineskip}\nopagebreak
\begin{center}
\large \bf 2.3. Connection of Hyperbolic Inflation with Lorentz Transformations of Minkowski Space
\end{center}
\vspace{.2cm}
\rm Using Lorentz transformations we give an alternative definition of hyperbolic expansion:\\
\\
\rm As is known, Minkowski space is a pseudo-Euclidean space of signature $(1,3)$, where the squared interval between $(t,x,y,z)$ and $(t_1, x_1, y_1, z_1)$ is given by $c^2(t-t_1)^2 - (x-x_1)^2-(y-y_1)^2-(z-z_1)^2$. Set $c = 1$. The Lorentz transformation of a frame moving with speed $v<1$ along the $x$-axis is $(x,y,z,t) \rightarrow \left( \frac{x-vt}{\sqrt{1-v^2}}, y,z,\frac{t-vx}{\sqrt{1-v^2}}\right)$. The squared interval is invariant under Lorentz transformations. Restrict the Lorentz transformations and the function $(t-t_1)^2 - (x-x_1)^2-(y-y_1)^2-(z-z_1)^2$ to the hyperplane $z=0$. After restriction one can identify points of the resulting space with points of the pseudo-Euclidean space of signature $(2,1)$, which in turn correspond to oriented circles. Thus Lorentz transformations induce transformations of the set of cycles: $\mathscr{L}_v :(x,y,r) \rightarrow \left( \frac{x-vr}{\sqrt{1-v^2}}, y,\frac{r-vx}{\sqrt{1-v^2}}\right)$. \\
\rm Study the properties of this transformation: if the center of a cycle $\omega$ lies on a line $y=\text{const}$, then the center of $\mathscr{L}_v(\omega)$ also lies on it. We prove that $x=0$ is the radical axis of these circles. Indeed, the point $(0,0)$ has the same power with respect to them since $x^2+y^2-r^2 = \left(\frac{x-vr}{\sqrt{1-v^2}}\right)^2+y^2-\left(\frac{r-vx}{\sqrt{1-v^2}}\right)^2$. Since $x=0$ is perpendicular to the line of centers, it is the radical axis. Therefore, if a circle lies entirely in the half-plane $x>0$, then in the Poincaré model in this half-plane $\omega$ and $\mathscr{L}_v(\omega)$ are concentric. To prove that their radii differ by a constant depending only on $v$, it suffices to compute the difference of the lengths of the diameters cut by the perpendicular from the cycle center to $x=0$. The intersections with the cycle corresponding to $(x,y,r)$ are $(x-\sqrt{r},y)$ and $(x+\sqrt{r},y)$; the intersections with $\left( \frac{x-vr}{\sqrt{1-v^2}}, y,\frac{r-vx}{\sqrt{1-v^2}}\right)$ are $\left( \frac{x-vr}{\sqrt{1-v^2}} - \sqrt{\frac{r-vx}{\sqrt{1-v^2}}},y\right)$ and $\left( \frac{x-vr}{\sqrt{1-v^2}} + \sqrt{\frac{r-vx}{\sqrt{1-v^2}}},y\right)$. The hyperbolic diameter of the first is $\ln \left( \frac{x+r}{x-r}\right)$. The diameter of the second is:
$$
\ln\left(\frac{ \frac{x-vr}{\sqrt{1-v^2}} + {\frac{r-vx}{\sqrt{1-v^2}}}}{ \frac{x-vr}{\sqrt{1-v^2}} - {\frac{r-vx}{\sqrt{1-v^2}}}}\right) = \ln\left(\frac{1-v}{v+1}\right)+\ln\left(\frac{x+r}{x-r}\right)
$$
\\
As we see, the radius in hyperbolic geometry changed by a constant. Thus on the set of circles lying in the half-plane $x>0$ or $x<0$, the transformation $\mathscr{L}_v$ acts as inflation of circles in hyperbolic geometry. A similar computation verifies the coincidence of the transformations on the set of equidistant curves and horocycles. Therefore, after composing with inversion in a point not lying on $x=0$, this transformation becomes hyperbolic inflation. Note that this definition is much simpler to use than the previous one. In particular, the fact that tangent cycles map to tangent cycles follows from invariance of the squared interval under Lorentz transformations. Thus, if the squared interval is $0$, corresponding to tangent circles, then after the transformation it remains $0$.\\
\\
Comment:
The authors did not manage to identify the physical nature of the relationship between inflation of cycles in hyperbolic geometry and Lorentz transformations, although its existence may well correspond to some physical phenomenon.\\
The transformations $\mathscr{L}_v$ were introduced in [1], albeit using exclusively Euclidean techniques.\\
\\
\Needspace{6\baselineskip}\nopagebreak
\begin{center}
\large \bf 2.4. Proof of the Main Theorem\\
\end{center}
\vspace{.2cm}
\rm We use the following lemma from [4] (Lemma 1): if there are three conics $\mathscr{C}_1, \mathscr{C}_2,\mathscr{C}_3$ not lying in one pencil and conics $\mathscr{F}_{12},\mathscr{F}_{23}$ such that $\mathscr{F}_{12}$ lies in the pencil of $\mathscr{C}_1,\mathscr{C}_2$, and $\mathscr{F}_{23}$ lies in the pencil of $\mathscr{C}_2,\mathscr{C}_3$, then there exists a conic $\mathscr{F}_{13}$ lying simultaneously in the pencil of $\mathscr{C}_1, \mathscr{C}_3$ and of $\mathscr{F}_{12},\mathscr{F}_{23}$.\\
\rm Using this lemma we show that any circle $\Omega$ tangent to $\gamma$ at two points lies in one pencil with $(ABCD)$ and with some circle tangent to $AB$ and $CD$. Indeed, let the tangency points be $X,Y$, and set $\mathscr{C}_1, \mathscr{C}_2,\mathscr{C}_3$ to $(ABCD)$, $\gamma$, $\Omega$. Let $\mathscr{F}_{12},\mathscr{F}_{23}$ be the product of the lines $AB, CD$ and the square of the line $XY$. By the lemma there exists a conic lying in the pencil with $(ABCD)$ and $\Omega$, and also in the pencil of $AB, CD$ and the square of $XY$. Since this conic lies in a pencil with two circles, it is a circle. Hence there exists a circle tangent to the lines $AB, CD$ at their intersection points with $XY$ and lying in the pencil with $(ABCD)$ and $\Omega$. 
\\
\rm From the argument above it follows immediately that the axes of symmetry of the conic $\gamma$ are parallel to the angle bisectors of the angle formed by $AB, CD$, since both pairs are perpendicular/parallel to $XY$. Consider the locus of centers of circles tangent to $\omega, \omega_1$ with one internal and the other external tangency. This locus is an ellipse with foci at the centers of the circles. Its interior clearly contains the intersection of the disks bounded by $\omega, \omega_1$, which contains the center of $\gamma$. Therefore the axes of symmetry intersect this locus, and hence there exists a circle tangent to $\gamma$ twice whose center lies on this locus. Thus we have obtained two concentric circles such that for both there exists a circle from the pencil formed by it with $(ABCD)$ that is tangent to $AB$ and $CD$. It is easy to check that the centers of the circles found lie on the same angle bisector of $AB$ and $CD$. But then the circles coincide, since their centers are determined as the intersection of the line joining the center of the obtained circles with the center of $(ABCD)$. Hence the circles obtained coincide, which yields the main theorem.
\\
\\
\Needspace{6\baselineskip}\nopagebreak
\begin{center}
\large \bf 3. Acknowledgements
\end{center}
\it The authors would like to thank Aleksey Suvorov and Roman Prozorov for their valuable comments on the work.\\
\Needspace{6\baselineskip}\nopagebreak
\begin{center}
\large \bf 4. Bibliography
\end{center}
\it [1] Isaak Yaglom, ``Geometric Transformations'', Part II, State Publishing House of Technical Literature, 1956\\
\it [2] Aleksey Suvorov, ``Three-Dimensional Space with One Imaginary Coordinate and Casey’s Theorem'', MMKSh, 2021\\
\it [3] P.V. Bibikov, I.I. Frolov, ``Non-Euclidean Solutions for Euclidean Problems'', Mathematical Enlightenment, Issue 26\\
\it [4] V.D. Konyshev, ``On Conics Tangent Twice and on Desargues’s Involution Theorem'', Mathematical Enlightenment, Issue 35\\
\it [5] Belavin A.A., Kulakov A.G., Usmanov R.A., ``Lectures on Theoretical Physics'', 2nd ed., MCCME, 2001\\
\end{document}